\documentclass{article}

\usepackage{latexsym}
\usepackage{amssymb}
\usepackage{amsmath}

\newtheorem{theorem}{Theorem}[section]

\newtheorem{corollary}[theorem]{Corollary}
\newtheorem{proposition}[theorem]{Proposition}
 \newtheorem{thmsn}{Main Theorem}

\newtheorem{definition}[theorem]{Definition}

\numberwithin{equation}{section}

\def\D{{\mathcal D}}
\def\E{{\mathcal E}}

\def\R{{\mathbb R}}
\def\V{{\mathcal V}}

\def\P{{\mathcal P}}
\def\U{{\mathcal U}}

\def\L{{\mathcal L}}

\def\hR{{}^{*}\R }
\def\X{{}^{*}\! X }

\def\A{{}^{*}\! A }
\def\B{{}^{*}\! B }
\def\f{{}^{*}\! f }
\def\g{{}^{*}\! g }

\def\p{{}^{*}\! p}
\def\q{{}^{*}\! q}
\def\ns{{}^{*}\!}

\def\ac{\check{\alpha}}
\def\ag{\alpha}
\def\ah{\widehat{\alpha}}
\def\bh{\widehat{\beta}}

\def\fs{{\mathsf f}}

\def\U{{\mathcal U}}

\def\Ux{{\mathcal U}_{\xi}}
\def\V{{\mathcal V}}

\def\0{\emptyset}
\def\/{\setminus}
\def\_{\overline}
\def\*{\times}
\def\+#1{\vec{#1}}

\def\ult#1#2{^{#1}_{\; #2}}
\def\lult#1#2#3{^{#1}_{\; #2}|#3}
\def\incl{\subseteq}

\def\Iff{\Longleftrightarrow}

\def\pes{\emph{e.g.}}

\def\ie{\emph{i.e.}}
\def\Ie{\emph{I.e.}}

\def\qed{\hfill $\Box$}

\begin{document}

\title{A simple algebraic characterization of nonstandard
extensions\thanks{Work
  partially supported  by  MIUR Grants PRIN 2007,2009, Italy.
  ${~~~~~~~~~~~~~~~~~~~~~~~~~~~~~~}$
  2010 MSC: Primary  03H05; 03C07; 03C20;
 Secondary: 26E35.
 }}
\author{Marco Forti\\
Dipart. di Matematica Applicata ``U. Dini'' \\
Via Buonarroti 
1C - 56100 PISA (Italy)\\
{\tt forti@dma.unipi.it}}

\date{}
\maketitle


%
%



\begin{abstract}
    We introduce the notion of \emph{functional extension} of a set $X$, by
    means of two natural algebraic properties of the  operator ``$*$'' on unary
    functions. We
    study the connections with ultrapowers of structures with 
    universe $X$, and 
    we give a simple  characterization of those functional
    extensions  that correspond
     to limit ultrapower extensions. In particular 
    we obtain a purely algebraic proof of Keisler's characterization of 
    nonstandard
    (= complete elementary) extensions.

\end{abstract}

\section*{Introduction}

Besides the ``superstructure approach'' proposed in \cite{RZ} 
(see Section 4.4 of \cite{CK}), various presentations of the 
``nonstandard methods'' that use \emph{nonstandard set-theories}
have been started in
\cite{Ne} and \cite{Hr}, and are continuously developing. However it
seems that, in some sense, the nonstandard methods do not really need
such set-theoretic generality. In fact,
several different ``elementary'' approches are also available, 
starting from
 the very interesting one of \cite{Ke76} (see also  \cite{BDNax}).
A general ``algebraic'' characterization of nonstandard extensions is
given by  W.~S.~Hatcher in \cite{Hat85}.
 Another purely algebraic construction producing every nonstandard
model is proposed in \cite{BDNal}.
A survey \cite{8path} of such ``elementary introductions'' to  nonstandard
methods has been jointly presented by 
V.~Benci, M.~Di Nasso, and the author.

This paper originates from a reflexion about the topological approach 
to nonstandard
models  exploited in \cite{brasil,DNFtop}. There a
nonstandard extension of a set $X$ is considered 
as a sort of
``topological completion'' $\,\X$, where $X$ is a discrete dense
subspace and each  function $f:X\to X$ has a continuous
extension $\f:\X\to \X$. In the  Hausdorff case of \cite{brasil}, 
these extensions are uniquely determined by the topology.
Since these nonstandard models turned out to be
a very narrow class of spaces,  more
general $T_{1}$-spaces are considered in \cite{DNFtop}. 
Uniqueness of continuous extensions of
functions being so lost, the
    ``$\ast $'' operator was charged to provide a
\emph{distinguished continuous extension} to $\X$ of each  function $f:X\to X$.
As a consequence,  in this general case the topology does not
force the choice of the continuous extensions, but it is rather this choice
that induces a topology on $\,\X$. These considerations suggest that one
might find
purely ``algebraic'' conditions on the $\ast$-extensions of functions, 
so as to characterize  all
nonstandard extensions, without any mention of
topologies.

In this paper we consider simple supersets $\,\X$ of $X$ together
with
an operator $*:X^{X}\to \X^{\X}$, which provides a
\emph{distinguished  extension} of each  function $f:X\to X$. We
  show that three algebraic conditions on the ``$\ast $''
operator are all what is needed to make such a functional extension
a true \emph{nonstandard model} of $X$,
and this in a very simple and natural way. Our main result is the
following
\begin{thmsn}
Let $\,\X$ be a  superset of the set $X$, and assume that to
every function $f:X \to X$ is associated a distinguished
extension $\f:\X \to \X$ satisfying the following conditions:

\begin{itemize}
\item[$(\mathsf{comp})$]  $~\g\circ \f = \ns(g\circ f),$ ~for all $f,g:X
\to X$;

\smallskip
	 \item[$(\mathsf{diag})$]~ $\ns(\chi\circ(f,g))(\xi)= \left \{ 
\begin{array}{ll}
			     1 & \mbox{\emph{if~}~} \f(\xi)=\g(\xi)\\
			      0 & \mbox{\emph{otherwise}}
			       \end{array}
		     \right.$  \ \ \ \ for all $f,g:X \to X$ and
all $\xi\in\X$, where $\, \chi: X\times X \to \{0,1\}$ is the characteristic
function of the diagonal, i.e. $\chi(x,y)=1 \iff x=y$.

\smallskip
\item[$(\mathsf{dir})$]~  for all $\xi,\eta \in \X$ there
exist $f,g:X\to X$ and $\zeta \in \X$ such that
$\ \f(\zeta) = \xi\,$ and $\,\g(\zeta) = \eta$;
\end{itemize}

\smallskip
Then there exist an ultrafilter $\D$ over $I=\X\* X$, and a filter 
$\E$ of equivalences on $I$ such that $\X$ is isomorphic to the limit 
ultrapower $X\lult{I}{\D}{\E}$. So $\X$ is a nonstandard extension of 
$X$.

\end{thmsn}

For precise definitions of the involved notions we refer to the next 
section. Here we simply anticipate that, like it is apparent for 
$(\mathsf{comp})$ and $(\mathsf{diag})$, also the third condition
$(\mathsf{dir})$ admits a (stronger) \emph{first order formulation}. So, on the one 
hand our criterion generalizes and greatly simplifies that of 
\cite{Hat85}. On the other hand, by directly defining the ultrafilter $\D$ 
and the equivalences $\E$, it provides a new, ``purely algebraic'' proof of 
the celebrated Keisler's characterization of \emph{complete elementary 
extensions} (see Theorem 6.4.10 of \cite{CK}).
We hope that
this axiomatic approach to nonstandard methods may be of independent interest.

\smallskip
The paper is organized as follows.
In Section \ref{uno}, we give the formal definition of functional
extensions and we study the first ``preservation properties'' of these
extensions.
In Section \ref{due}, we 
prove the main theorem and obtain our ``functional'' characterization of nonstandard 
($=\,$complete elementary)
extensions.
Concluding remarks can be found in the final Section \ref{froq}, 
where in 
particular we
briefly outline the connections between
functional extensions and the topological extensions of
\cite{DNFtop}.

In general, we refer to
  \cite{CK} for definitions and facts
concerning ultrapowers, ultrafilters, and nonstandard models that are 
used in this paper.

\medskip
 The author is grateful to  Mauro Di Nasso and Vieri Benci for useful 
 discussions and suggestions.

\section{Functional extensions}
\label{uno}

A main feature of all nonstandard models of Analysis is the
existence of a canonical extension $\f:\hR \to \hR$  of any (standard) function
$\, f:\R \to \R$ (and also of any subset $A\incl \R$). Here we use this
property as the definition
of \emph{functional extensions} of an arbitrary set $X$.
  We shall assume in the sequel that
$0,1 \in X$, in order to have at disposal the extensions of 
\emph{characteristic
functions}.  Since the extension
of the characteristic function of $A\incl X$
  will turn out to be a characteristic function in $\,\X$, we shall use it to
  define  the \emph{extension} $\A$ of $A$ in $\,\X$.

The so called \emph{nonstandard methods} are intended to study
extensions which preserve those properties of the standard structure
which are currently being considered.
The \emph{Transfer $($Leibniz's$)$ Principle} states that all properties
that are expressible
in a(\emph{n in})sufficiently expressive language are preserved by passing to
the nonstandard models.
In particular, an important
property of nonstandard models of Analysis is that
``disjoint functions have disjoint extensions''.
This property has a clear ``analytic'' flavour, and in fact it can be
considered as the most
characteristic
feature of \emph{nonstandard} extensions when
compared with \emph{continuous} extensions of functions in
compactifications or topological completions,
where  equality may be reached  at limit points only (see 
\cite{DNFtop}).
So we should reasonably postulate  this property, together with
preservation of characteristic functions.

With this in mind, in the following definition we assume that
``composition and diagonal are
preserved'':

\begin{definition}\label{fext}
A  superset  $\,\X$ of the set $X$ is a
\emph{functional extension } of $X$ if to
every function $f:X \to X$ is associated a distinguished
\emph{extension } $\f:\X \to \X$ in such a way that the following 
conditions are fulfilled, 
 for all
$f,g:X \to X$ and all
$\xi\in\X:$

\begin{itemize}
\item[$(\mathsf{comp})$]  $~~\g(\f(\xi))= \ns(g\circ f)(\xi),$~
\end{itemize}

and

\begin{itemize}
\item[$(\mathsf{diag})$]~ $\ns(\chi\circ(f,g))(\xi)= \left \{ \begin{array}{ll}
			     1 & \mbox{\emph{if~}~} \f(\xi)=\g(\xi)\\
			      0 & \mbox{\emph{otherwise}}
			       \end{array}
		     \right. $
\ \ \ where $\, \chi: X\times X \to \{0,1\}$\\ is the characteristic
function of the diagonal, i.e. $\chi(x,y)=1 \iff x=y$.

\end{itemize}

\smallskip\noindent
The functional extension $\X$ of $X$   is \emph{directed} if 
\begin{itemize}
\item[$(\mathsf{dir})$]~  for all $\xi,\eta \in \X$ there
exist $f,g:X\to X$ and $\zeta \in \X$ such that
$\ \f(\zeta) = \xi\ $ and $\g(\zeta) = \eta$.
\end{itemize}

\end{definition}

  Following the common usage, we call
\emph{standard} the points of $X$
and \emph{nonstandard} those of $\,\X \/ X$. The adjective 
``directed'' refers to the so called \emph{Puritz order}, 
a preordering of nonstandard models corresponding to the 
Rudin-Keisler (pre)ordering of ultrafilters (see \cite{pu,ng}): 

\begin{center}
    for $\xi,\eta \in \X$ put $\eta\le_{P}\xi$ if there exists 
    $f:X\to X$ such that $\eta=\f(\xi)$.
\end{center}
Then clearly the property 
$(\mathsf{dir})$ expresses that the Puritz order $\le_{P}$ on $\X$
is \emph{directed}.

\smallskip
We are interested only in directed extensions, but all functional extensions 
satisfy various natural ``preservation
properties''. We begin by considering characteristic functions, so as
to define the extensions of subsets and derive their properties.

\begin{theorem}\label{bool} ~Let $\,\X$ be a functional extension of $X$.
Then
\begin{itemize}
\item[$(i)$] if $\, c_{x}: X \to X$ is the constant function
with value $x\in X$, then the extension
$\, \ns {c_{x}}$ is the constant with value $x$ on $\, \X$;
\item[$(ii)$] if $\, \chi_{A}: X \to X$ is the characteristic
function of  $A\incl X$, then the  extension
       $\ns {\chi_{A}}$ is the characteristic
function of a superset $\A$ of $A$ in $\, \X$.

\item[$(iii)$]
The  map ~$\ns :A\mapsto \A$ commutes with binary union, intersection and
complement. Moreover $\,\A\cap X = A$, hence $\ns$ is a boolean isomorphism
of  $\P (X)$ onto a subfield $\,St(\X)$ of  $\,\P (\X)$.

\end{itemize}
\end{theorem}

{\bf Proof.}~
Since $c_{1} = \chi\circ(f,f)$, the property $(\mathsf{diag})$ implies
that $\ns c_{1}$ is the constant $1$ on $\,\X$. Any other constant is
obtained from $c_{1}$ by composing with a suitable transposition,
and so $(i)$ follows from $(\mathsf{comp})$.

  Moreover
$\chi_{A} =  \chi\circ(\chi_{A},c_{1})$, hence   $\ns \chi_{A}$ takes
only the values $0$ and $1$. Since obviously it maps $A$ to $1$, we
have $(ii)$.

  Let $\tau$ be a transposition mapping
$0$ to some $x\ne 0,1$. Then $\ns (\tau\circ\chi_{B})$ is $1$ exactly
on $\B$, and so
$\ns (\chi_{A\cap B}) = \ns (\chi\circ(\chi_{A},\tau\circ\chi_{B})) =
\chi_{\A\cap\B}$.
So  the  map $\ns$ commutes with binary intersection. Moreover 
$\chi\circ(\chi_{A},\chi_{X\/ A}) = c_{0}$, hence $\A$
and $\ns (X\/ A)$ give a partition of $\X$. So we have at once that
$\,\A\cap X = A$ and that $\X\/ \A= \ns(X\/ A)$. It follows that the 
map $\ns$
is a boolean isomorphism, and $(iii)$ is proved.

\hfill $\Box$

\bigskip
Having defined $*$-extensions of sets, the property $(\mathsf{diag})$
gives immediately  a sort of ``preservation of
equalizers'', which corresponds to another basic idea of nonstandard
analysis, namely that
``standard functions behave like germs'':

\begin{corollary}\label{preseq}
Let $\, \X$ be a functional extension   of
$X$. Then
$$\{ \xi \in \X \mid \f(\xi)= \g(\xi)\} =
\ns{\{ x \in X \mid f(x)= g(x)\}}$$
for all $f,g: X \to X$; or equivalently, for all $\xi\in\X$,
$$\f(\xi)= \g(\xi)\ \Longleftrightarrow \
	\exists A\incl X. \,( \xi \in \A \ \ \& \ \forall x\in A \ .
f(x)=g(x)\,).$$
\hfill $\Box$
\end{corollary}

\medskip
It follows from  Theorem \ref{bool} and Corollary \ref{preseq}  that 
$\X$ induces a
uniquely determined functional extension $\, \A$  of any subset
$A\incl X$. Namely,  for $f:A \to A$,  let $\f:\A \to \A$ be
the restriction to $\A$ of  $\, \g$, where $g:X\to X$ is any function
whose restriction to  $A$ is $f$.

Notice that the identity map \emph{may not be preserved} by
functional extensions.  If
	  $\imath:X\to X$ is the identity of $X$, in the general case one
only obtains

\begin{center}
	$ \f (\xi)= \f (\ns \imath(\xi))= \ns \imath( \f
(\xi))~$ for all  $f:X\to X$
and all $\xi
	\in \X$.
\end{center}

\noindent Thus $\ns\imath$ is the identity exactly on those
points of $\X$ which are  reached by some function $\f$, and all
functions $\f$ map each nonstandard point $\xi$ to the same point as
$\ns\imath(\xi)$. So, when $\ns\imath$ is not the identity,
the extension $\X$ can be considered ``redundant'', in the sense that
the extensions of all functions are completely determined by
their restrictions to $\ns \imath(\X)$, and the remaining elements of
$\,\X$ are not attained by any function $\f$. Moreover $\ns \imath(\X)$,
equipped with the restrictions of all  $\f$s, becomes a functional
extension where the identity is preserved.

Call \emph{irredundant} a functional extension $\,\X$ of $X$ if
the following weakening of $\mathsf{(dir)}$ holds:
\begin{center}
for all $\xi\in\X$ there exist $f:X\to X$ and
$\eta\in\X$ such that $\f(\eta) = \xi.$
\end{center}
Important and natural  preservation properties, concerning
\emph{ranges,} \emph{injectivity,} and \emph{finite subsets} can be derived
for irredundant extensions, hence \emph{a fortiori} for directed extensions:

\begin{proposition}\label{presirr}
Let $\, \X$ be an irredundant functional extension   of
$X$. Then, for all $f: X \to X$ and all $A\incl X$
\begin{itemize}

\item[$(i)$] $\, \f(\A) = \ns(f(A))\,$
    (in particular $\f$ is
onto if $f$ is onto);
\item[$(ii)$] ~if $f:X\to X$ is one-one on $A$,  then
    $\, \f$ is
one-one on $\A$.

\item[$(iii)$] Extensions of finite sets are trivial, i.e.
	$A=\A$ whenever $A$ is finite.

\end{itemize}

\end{proposition}

{\bf Proof.}~~
In irredundant extensions $\ns\imath$ is the identity. Hence $(i)$ and 
$(ii)$ follow from
$(\mathsf{comp})$, because (the restriction of) a function is injective (resp.
surjective) if and only
if it has a left (resp. right) inverse.

In order to obtain $(iii)$ observe first that $\ns\chi_{\{1\}} =
  \ns(\chi\circ(\imath ,c_{1}))$ is $1$ exactly on $1$. Hence $\ns\{1\}=
  \{1\}$, and the same property holds for all singletons, by point
  $(i)$. Since $\ns\,$ is a boolean isomorphism, all of $(iii)$ follows.

\hfill $\Box$

\medskip
In order to make an effective use of  nonstandard models,
it is always  assumed by nonstandard analysts that \emph{all and only 
infinite sets}
are
indeed \emph{extended}, i.e. $A=\A$ if and only if $A$ is finite. We
shall not need this assumption, but all properties of Proposition
\ref{presirr} are valid in any nonstandard model.
So we have in mind essentially  only irredundant extensions.  We have not
postulated
irredundancy in Definition \ref{fext} because
this condition is still too weak to obtain \emph{fully nonstandard} 
extensions. To this aim we have isolated in the Introduction the stronger
property $\mathsf{(dir)}$ that
``every \emph{pair}
of points is dominated by some  point''.

Remark that all the properties stated in Theorem \ref{bool} and its corollaries, as well as the
defining properties $(\mathsf{comp})$ and
$(\mathsf{diag})$, are particular
(and
\emph{perspicuous}) cases of the Transfer Principle. On the other 
hand, this seems \emph{prima
facie} not to apply
to the condition $(\mathsf{dir})$, whose straightforward formalization
is \emph{second-order}.  On the contrary, a
\emph{strong uniform version} of directedness can be obtained as an
instance of
transfer. Namely one can compose any bijective function  
$\delta :X\to X\times X \ $
with the ordinary projections $\pi_{i}: X\times X \to X$, so as to 
obtain
``unary projections''
$p_{1},p_{2}:X\to X$ such that
\begin{center}
 \emph{for all $x,y\in X$ there is a
unique $z\in X$ such that
   $p_{1}(z) =x,\ p_{2}(z) =y$.}
\end{center}
Then by Transfer one obtains the following strengthening of 
$(\mathsf{dir})$:
\begin{center}
 \emph{for all $\xi,\eta\in\X$ there is a
unique $\zeta\in\X$ such that
   $\p_{1}(\zeta) =\xi,\ \p_{2}(\zeta) =\eta$.}
\end{center}

So, if we want a full Transfer Principle in our functional
extensions,
then  not only irredundancy, but also
\emph{directedness} has to hold at the end. On the
other hand, the very point of this paper is the remarkable fact that 
we  shall prove in the next section:\\
$\bullet$ \emph{by combining the \emph{sole}
property $(\mathsf{dir})$ with the \emph{simple,
natural} conditions $(\mathsf{comp})$ and $(\mathsf{diag})$,
one forces the
strongest \emph{Transfer Principle for all first order properties},
thus providing  \emph{complete elementary extensions}.
}

\section{Algebraic characterization of directed functional extensions}
\label{due}

We devote this section to a ``purely algebraic'' proof of the Main 
Theorem.

Recall that two functions $f,g:I\to X$ are equivalent modulo $\,\U$, where
$\U$ is an arbitrary ultrafilter over $I$,
if they agree on some set $U\in\U$. The ultrapower $X\ult{I}{\U}$
is the set of the
equivalence classes modulo $\,\U$ of all functions $f:I\to X$. We
refer to \cite{CK} for basic facts about ultrapowers. In the sequel,
in dealing with ultrapowers, we shall adhere to the following notation:
      \begin{itemize}
	\item $[f] \in X\ult{I}{\U}$ is the equivalence class of the
	function
      $f:I\to X$;
	\item  $\_g: X\ult{I}{\U}\to  X\ult{I}{\U}$ is
      the \emph{interpretation} of the function $g:X\to X$ in the 
      ultrapower, \ie\
      $\_g ([f]) = [g\circ f]$ for all $f:I\to X$;
	\item  $ A^{I}/\U \incl X\ult{I}{\U}$ is the \emph{interpretation} of
$A\incl X$ in the
	ultrapower.
	\end{itemize}

\smallskip
The {\it limit ultrapower}
$X\ult{I}{\U}|\E$, where $\U$  is an ultrafilter over $I$ and
$\E$ is a filter of equivalences on $I$, is the
subset of the ultrapower $X\ult{I}{\U}$ containing
the $\U$-equivalence classes of all functions
$f:I\to X$ that induce on $I$ an equivalence $$Eq(f) = \{(x,y) \mid
f(x)=f(y) \}$$ that belongs to $\E$.
(Notice that one can assume
w.l.o.g.~that
each $E\in\E$ gives a partition of size not
exceeding $|X|$.)

A celebrated  theorem of Keisler's  states that  every \emph{complete 
elementary extension}\footnote{~\Ie\ an extension that satisfies the 
Transfer Principle w.r.t. a language containing \emph{all $n$-ary relations} 
over $X$.}
of a structure $\mathfrak
X$ with universe $X$ is isomorphic to a  limit ultrapower
  of $X$ (see \pes\ Section 6.4 of \cite{CK}). We shall obtain 
  the same conclusion from the much weaker assumption of the three 
  properties $(\mathsf{dir})$, $(\mathsf{comp})$, and $(\mathsf{diag})$.
  
 Starting from a given directed functional extension
$\,\X$ of $X$, we now
define a limit ultrapower
$X\lult{I}{\D}{\E}$ that will result isomorphic to $\X$.

\begin{itemize}
    \item For $\alpha\in\X$ put
$\check{\alpha} = \{\xi \in \X \mid \exists f:X\to X\ s.t.\ 
\alpha=\f(\xi)\}$, and\\
for all $\xi \in \check{\alpha}$ fix a
function $f_{\xi\alpha}:X\to X$ such that $\f_{\xi\alpha}(\xi) = \alpha$.

\smallskip
\item The family $\{\ac \mid \alpha \in 
\X \}$ has
the \emph{finite intersection
property}, because $\,\X$ is directed, and so we can extend it to an ultrafilter
   over $\,\X$, say $\V$.

   \smallskip
\item Put $\ I = \X \times X$.  For $\alpha \in \X$,  define
the function 

  $$\ah : I\to X \ \ \ \mbox{by}\ \ \ \  \ah (\xi,x) = \left \{ \begin{array}{ll}
	 f_{\xi \alpha}(x) & \mbox{if~~} \xi \in \ac\\
			     x & \mbox{otherwise}
			       \end{array}
		     \right.$$
and the corresponding equivalence $E_{\alpha}$ induced  on $I$ by 
$$E_{\alpha} = Eq(\ah)=\{(i,j)\in I\times I \mid \ah (i) =\ah(j)\}.$$

\smallskip
\item Then $E_{\alpha}\incl E_{\beta} \cap E_{\gamma}$ whenever 
$\ag \in \check{\beta}\cap\check{\gamma}$, and so the family 
 $\{E_{\alpha}\mid \ag\in\X\}$ generates a filter of equivalences
on $I$, say $\E$.

\smallskip
\item Let $\, \Ux=\{A\incl X\mid \,\xi\in\A\}$ be the
ultrafilter on $X$ generated by $\xi$,\hfill

\smallskip
\noindent and let $\D =
\sum_{\V}\,\U_{\xi}$ be the ultrafilter on $I$ such that
$$D\in\D \ \iff \ \{\xi\in\X \mid \{x\in X \mid (\xi,x)\in D \}\in
\U_{\xi} \}\in \V.$$

\end{itemize}

Then we have


\begin{theorem}\label{funalg}~
Let $\,\X$ be a  directed functional extension of $X$. Then
the map $$\widehat{}\,:\X \to X^{I}$$ induces an isomorphism of $\,\X$ onto the 
limit ultrapower $X\lult{I}{\D}{\E}$ (as
structures for the language $\L = \{f \mid f: X\to{X}\}$).
\end{theorem}

{\bf Proof.}~
First of all we claim that, for all  $\alpha, \beta \in
\X$ and all $g:X\to X$, $$\beta = \g(\alpha)\ \ \Iff\ \ 
\ [\widehat{\beta}\,] = [g\circ \widehat{\alpha}]\ \ (\!\!\!\!\!\mod \D).$$

Assume $\beta = \g(\alpha)$. Then, for all $\xi\in\ac$, $\,\f_{\xi
\beta}(\xi) = (\g\circ\f_{\xi \alpha})(\xi)$. Applying  Corollary
\ref{preseq} we obtain
  $$\, \{x\in X \mid (g\circ f_{\xi\alpha})(x) = f_{\xi\beta}(x)
\} = \{x\in X \mid  g(\ah (\xi,x)) = \bh (\xi,x) \}
\in \Ux$$ for all $\xi\in\ac$. Hence $\{i\in I \mid g\circ\ah = \bh
\} \in \D$, \ie\ $\ [\widehat{\beta}\,] = [g\circ \widehat{\alpha}]$.

\smallskip
Conversely, assume $\{i\in I \mid g\circ\ah = \bh \} \in \D$. Then

\smallskip
$~~~~~~~~~\{\xi\in\ac \mid \{x \mid  g(\ah (\xi,x)) = \bh (\xi,x) \}\in
\Ux\} =$

\smallskip
$~~~~~~~~~ \{\xi\in\ac \mid [f_{\xi\beta}]=[g\circ f_{\xi\alpha}]\ \mathrm{mod}
\ \Ux \} =~~~\{\xi\in\ac \mid \f_{\xi\beta}(\xi) =
\g(\f_{\xi\alpha}(\xi) \} \in \V.$

\smallskip
\noindent Therefore $\beta = \g(\alpha)$.

\smallskip
It follows at once that the map $\ \widehat{} \ $ is injective modulo 
$\D$, and that it
preserves
extensions of functions. We prove now that the range of 
$\ \widehat{} \ $ contains all functions 
$\varphi:I\to X$ such that  $Eq(\varphi)\in\E$.

\smallskip
Let $\varphi:I\to X$ be such a function. Then there
exists $\alpha \in \X$ such that $E_{\alpha}\incl Eq(\varphi)$, i.e.
$\varphi$ is constant on each
equivalence class
modulo $E_{\alpha}$. Therefore $\varphi$ is completely determined by
its behaviour on a set of representatives, which may be conveniently
chosen as $\{(\alpha,x) \mid x\in X\}$. 

Put $g(x) = \varphi
(\alpha,x)$ and $\beta = \g (\alpha)$. By the above claim we have
$$\widehat{\beta}(\alpha,x) = \widehat{\g(\alpha)} (\alpha,x) =
(g\circ \widehat{\alpha})(\alpha,x) = g(x).$$ Hence $\varphi =
\widehat{\beta}$.

\hfill $\Box$

\bigskip

According to the above theorem, any directed functional extension is a 
limit ultrapower, that in turn is a
complete elementary
extension. So we have the following corollary, which incorporates 
Keisler's Theorem:

\begin{corollary}\label{KT}
 The following conditions are equivalent for any
functional extension  $\, \X$   of
$X$:
\begin{enumerate}
    \item $\, \X$ is directed; 

    \item  $\, \X$ is a limit ultrapower of $X$;

    \item  $\, \X$ can be uniquely
    expanded to a nonstandard (complete elementary) extension of $X$.
\end{enumerate}
\qed
\end{corollary}

In the next section we shall sketch two different proofs, a 
``topological'' and a ``purely logical'' one,
of the implication $\,(1) \Rightarrow(3)\,$ in the above corollary, 
\ie\ the ``logical part'' of our Main Theorem. To be sure, 
 in order to completely 
prove Corollary \ref{KT}, either of these proofs has to be 
combined with Keisler's Theorem.

\section{Concluding remarks}
\label{froq}

We conclude this paper by outlining two different approaches, a 
``topological'' and a ``purely logical'' approach.  Each of them can 
provide a proof that directed 
functional extensions are complete elementary 
extensions.

\subsection{The Star topology of directed functional extensions}
\label{stop}
We pointed out in the Introduction that this paper has been inspired
by the topological approach to nonstandard models presented in
\cite{DNFtop}, where the extended functions $\f$ are taken
\emph{continuous} (and  $\A$ is the \emph{closure} in $\X$ 
 of $A\incl X$). However, while in the case of the 
 \emph{Hausdorff extensions}
of \cite{brasil},  the topology uniquely determines the extension, on 
the contrary, in the general
case it seems that the ``algebraic'' properties of the map $\ns$
are
responsible for the topology, rather than \emph{vice versa}.
Let us report the following definition from \cite{DNFtop}:

\begin{definition}\label{text}
A $T_{1}$ topological
space $\, \X$ is a
\emph{topological extension} of $X$ if

$\bullet$ $X$ is a \emph{dense} subspace of
$\, \X$, and 

$\bullet$ to
every function $f:X \to X$ is associated a distinguished
\emph{continuous extension} 

~$\f:\X \to \X$ in such a way that

\begin{itemize}
\item[$(\mathsf{c})$]  $  \g\circ \f = \ns(g\circ f)$~  {for all}
$f,g:X \to X$,
{and}

\item[$(\mathsf{i})$] if $f(x)=x$ for all $x\in A\incl X$, then
$\, \f (\xi)=\xi$ for all $\xi \in\_A$.
\end{itemize}

\noindent
$\bullet$ The topological extension $\, \X$ is \emph{analytic} if, for all
$f,g:X\to X$,

\vspace{1mm}
\begin{itemize}
\item[$(\mathsf{d})$]
~$f(x) \ne g(x)$ for all $x\in X$
$\ \ \Longrightarrow \ \ $
$\, \f(\xi) \ne \g(\xi)$ for all $\xi\in \X$.
\end{itemize}

\vspace{1mm}
\noindent
$\bullet$ The topological extension $\, \X$ is \emph{coherent} if
\begin{itemize}
\item [$(\fs)$] for all $\xi,\eta \in \X$ there exist
functions $p,q: X\to X$ and a point $\zeta \in \X$ such that
$\ \p(\zeta)=\xi$
and $\ \q(\zeta)=\eta.$
\end{itemize}

\end{definition}

It is easily seen that $(\mathsf{i})$ and $(\mathsf{d})$ follow from 
 $(\mathsf{diag})$. So 
all properties $(\mathsf{cidf})$ of the above
definition hold in any directed functional extension. As a matter of 
fact, it is
this definition that suggested the defining properties
$(\mathsf{comp,diag,dir})$ assumed in this paper.
So it seems natural to look for a topology that turns any 
\emph{directed functional
extension} into a \emph{coherent analytic
extension}.

In fact this task has been accomplished already in \cite{DNFtop}, where it is
shown how to put \emph{on every
nonstandard model} $\,\X$ a natural topology, named 
\emph{Star-topology}, that makes all functions $\f$
\emph{continuous}. The \emph{closed sets} of this topology are the
arbitrary intersections of  sets of the form
$$E(\+{f};\+{\eta}) = E(f_{1},\ldots,f_{n};\eta_{1},\ldots ,\eta_{n})
=\{ \xi \in \X \mid \exists i\in\{1,\ldots,n\}\,.\, \f_{i}(\xi) = \eta_{i},\,\} $$
for all $n$-tuples of functions $\ f_{i}:X\to X$, and of points $\eta_{i}
     \in \X$.

It is easily seen that the Star-topology is the coarsest $T_{1}$ topology
on $\,\X$ that makes all functions $\f$
continuous, and in fact we have

\begin{theorem}\label{funtop}~
Every directed functional extension  $\,\X$ of $X$, when endowed with the
Star-topology, becomes  a  coherent
analytic topological extension of $X$.

\end{theorem}

{\bf Proof.}~
We have already remarked that the properties $(\mathsf{cidf})$ are
 derivable from $(\mathsf{comp,diag,dir})$.
So we have only to prove that $X$ is dense in $\,\X$ w.r.t.~the
Star-topology.
Assume that $X\incl E(\+{f},\+{\eta})$. Since only nonstandard
elements can be mapped to nonstandard elements by functions $\f$, we
may assume w.l.o.g.~that $\eta_{i} =y_{i}\in X$ for $i=1,\ldots,n$.

Then $$E(\+{f},\+{y}) = \bigcup_{i}\f_{i}^{-1}(y_{i}) =
\bigcup_{i}\ns(f_{i}^{-1}(y_{i})) = \ns
(\bigcup_{i}f_{i}^{-1}(y_{i})).$$
It follows that $X = E(\+{f},\+{y})\cap X =
\bigcup_{i}f_{i}^{-1}(y_{i}),$ and so $\,\X=E(\+{f},\+{y})=\_X$.

\hfill $\Box$

\medskip

Now  Theorem 5.5  of
\cite{DNFtop} states  that all coherent analytic extensions are
complete elementary extensions. Thus, if we had been ready to accept
the full topological machinery of \cite{DNFtop}, we could
have obtained Theorem \ref{funalg} by combining the above
theorem with Keisler's Theorem. But the
aim of this paper is rather that of showing that \emph{few clear,
natural, purely
algebraic conditions} are all what is needed to obtain \emph{limit 
ultrapowers}, hence 
\emph{full nonstandard models}.

\subsection{An inductive logical proof}
\label{stop}

The very reason for considering only \emph{unary} functions in the
definition of functional extensions is the fact that the
properties $(\mathsf{comp,diag,acc})$ together allow for
a \emph{unique} \emph{unambiguous} definition of the
extension
$\, \ns\varphi$ of \emph{any $n$-ary function} $\varphi : X^{n} \to X$,
in a simple ``parametric'' way.

\begin{theorem}\label{extn}~ Let  $\, \X$ be an directed functional
extension of $X$.
Then there is a unique way of assigning an extension $\, \ns\varphi$ to
every $n$-ary
function $\, \varphi : X^{n} \to X$  so as to preserve all
compositions,\footnote{ ~i.e.~for all all $m,n\geq 1,$ all $\, 
\varphi : X^{n} \to
X$, and all
$\, \psi_{1}, \ldots, \psi_{n}  : X^{m} \to X$,

\vspace{-3mm}
$$\ns\varphi \circ (\ns\psi_{1}, \ldots,
\ns\psi_{n}) =  \ns(\varphi \circ (\psi_{1}, \ldots,
\psi_{n})).$$} namely
	$$\ns\varphi (\xi_{1},\ldots,\xi_{n})=
	\ns (\varphi \circ (f_{1},\ldots,f_{n}))(\zeta),$$
	where $f_{i}:X\to X$ and $\zeta\in\X$ are such that 
$\f_{i}(\zeta) = \xi_{i}$
	for $i=1,\ldots,n$.

\end{theorem}

{\bf Proof.}~ We can easily generalize the property $(\mathsf{dir})$ and prove
by induction on $n$ that, for all $\xi_{1},\ldots,\xi_{n}\in
\X$ there exist
	 $ f_{1},\ldots,f_{n}: X \to X$ and  $\zeta\in \X$ such that
$\f_{i}(\zeta) = \xi_{i}$
for $i=1,\ldots,n$.

When the  extensions of  $n$-ary functions
preserve compositions, 
the equality
$$\ns\varphi (\xi_{1},\ldots,\xi_{n})=
\ns (\varphi \circ (p_{1},\ldots,p_{n}))(\zeta)$$
has to hold whenever $\xi_{1},\ldots,\xi_{n}$, $ f_{1},\ldots,f_{n}$, and
$\zeta$ satisfy the above
conditions.
Therefore the extensions of \emph{unary} functions completely determine
those of all $n$-ary functions.

We are left with the task of showing that this
definition is independent of the choice of the functions $f_{i}$
and of the point
$\zeta$. Thus one has to prove
     $$\ns (\varphi \circ (f_{1},\ldots,f_{n}))(\xi) =
\ns (\varphi \circ (g_{1},\ldots, g_{n}))(\eta)$$
for all $\varphi : X^{n} \to X$, provided $
f_{1},\ldots,f_{n},g_{1},\ldots, g_{n}: X \to X$
and $\xi,\eta \in \X$ satisfy  $\f_{i}(\xi) = \g_{i}(\eta)$
for $i=1,\ldots,n$.

\smallskip
Taking a point $\zeta$ such that  $\p_{1}(\zeta)=\xi$ and
$\p_{2}(\zeta)=\eta$, we can assume
w.l.o.g.~that $\xi = \eta = \zeta$. Then,
applying Corollary \ref{preseq}, we get

\medskip
$~~\{ \xi\in\X \mid \ns (\varphi \circ (f_{1},\ldots,f_{n}))(\xi)
=  \ns (\varphi \circ (g_{1},\ldots, g_{n}))(\xi)\} = $

\medskip

$~~\ns{\{x\in X \mid (\varphi \circ (f_{1},\ldots,f_{n}))(x) =
(\varphi \circ (g_{1},\ldots, g_{n}))(x) \}} \supseteq $

\medskip

$~~\bigcap_{1\leq i\leq n}\ns{\{x\in X \mid  f_{i}(x) = g_{i}(x)
\}} = 
\bigcap_{1\leq i\leq n}\{\xi\in \X \mid  \f_{i}(\xi) = \g_{i}(\xi) \},$

\medskip
\noindent
and the proof is complete.

\hfill $\Box$

\medskip

\medskip
By using the characteristic functions in $n$ variables one can assign
an extension $\, \ns R$ also to each $n$-ary relation $R$ on $X$.
In this way, given a complete first order structure
$\frak X = \langle X; R_{i}\,,i\in I;F_{j}\,,j\in J \rangle$ and a
directed functional extension $\X$ of its universe, one produces
a similar structure
$\, \ns\mathfrak X = \langle \X; \ns R_{i}\,,i\in I;\ns F_{j}\,,j\in J 
\rangle$ with universe $\X$.

Now, in order to obtain a ``logical'' proof that 
$\mathfrak X$ is an \emph{elementary substructure} of $\, 
\ns\mathfrak X$, one should prove, for every formula $\Phi$, 
that $$\forall x_{1}, \ldots ,x_{n}\in X.\
	  {}^{*}\mathfrak X \models \Phi [x_{1}, \ldots ,x_{n}] \
	  \Longleftrightarrow \ \mathfrak X \models \Phi [x_{1}, \ldots
,x_{n}]. $$

The atomic case being dealt with in Theorem \ref{extn}, one
could proceed, as usual, by induction 
on the complexity of the formula $\Phi$. But developing the details 
of such an inductive proof lies outside the scope of this paper.

\bibliographystyle{amsplain}

\begin{thebibliography}{99}
    
\bibitem{BDNax}
\textsc{V.~Benci, M.~Di~Nasso} -
Alpha-Theory: an elementary axiomatics for nonstandard analysis,
\textit{Expo. Math.}, \textbf{21} (2003), 355--386.

\bibitem{BDNal}
\textsc{V.~Benci, M.~Di~Nasso} -
A purely algebraic characterization of the hyperreal numbers,
\textit{Proc. Amer. Math. Soc.} 133 (2005), 2501--2505.


\bibitem{brasil}
\textsc{V.~Benci, M.~Di~Nasso, M.~Forti} -
Hausdorff nonstandard extensions,
\textit{Boletim Soc. Parana.  Mat.} (3),
\textbf{20} (2002), 9--20.

\bibitem{8path}
\textsc{V.~Benci, M.~Di~Nasso, M.~Forti} -
The Eightfold Path to Nonstandard Analysis, in
    \emph{Nonstandard Methods and Applications in Mathematics}
    (N.J.~Cutland, M.~Di~Nasso, D.A.~Ross, eds.), L.N. in Logic \textbf{25},
    A.S.L. 2006, 3--44.


\bibitem{CK} \textsc{C.C.~Chang, H.J.~Keisler} -
\textit{Model Theory} (3rd edition),
North-Holland, Amsterdam 1990.


\bibitem{DNFtop}
\textsc{M.~Di~Nasso, M.~Forti} -
Topological and nonstandard extensions,
\textit{Monatshefte f\"ur Mathematik} \textbf{144} (2005), 89--112.

\bibitem{Hat85} {\sc W.S. Hatcher} - Elementary extension and the
hyperreal numbers,  in \textit{Mathematical Logic and
Formal Systems} L. N. Pure Appl. Math., Dekker, New York 1985, 205--219.

\bibitem{Hr} {\sc K. Hrb\`a\u{c}ek} - Axiomatic foundations for
Nonstandard Analysis, \emph{Fund. Math.} \textbf{98} (1978), 1--19.

\bibitem{Ke76} \textsc{ H.J.~Keisler} -
\textit{Foundations of Infinitesimal Calculus},
Prindle, Weber and Schmidt, Boston 1976.

\bibitem{Ne} \textsc{E.~Nelson} -
Internal Set Theory; a new approach to nonstandard arithmetic,
\emph{Proc. London. Math. Soc.} \textbf{83} (1977), 1165--1198.

\bibitem{ng} \textsc{S.-A.~Ng, H.~Render} -
The Puritz order and its relationship to the Rudin-Keisler order,
in \emph{Reuniting the antipodes - Constructive and nonstandard views of 
the continuum} (P.~Schuster, U.~Berger, H.~Osswald, eds.), Kluwer AP
Dordrecht 2001, 157--166.

\bibitem{pu} \textsc{C.~Puritz} -
Ultrafilters and standard functions in nonstandard analysis,
\emph{Bull. Amer. Math. Soc.} (3) \textbf{22} (1971), 705--733.

\bibitem{RZ} \textsc{A.~Robinson, E.~Zakon} - A set-theoretic
characterization of enlargements, in
\textit{Applications of Model Theory to Algebra, Analysis and Probability}
(W.A.J.~Luxemburg, ed.), New York 1969, 109--122

\end{thebibliography}

\end{document}